\documentclass[oneside,reqno]{amsart}
\usepackage[T1]{fontenc}
\usepackage[latin9]{inputenc}
\usepackage{color}
\usepackage[english]{babel}
\usepackage{prettyref}
\usepackage{mathrsfs}
\usepackage{enumitem}
\usepackage{amsbsy}
\usepackage{amstext}
\usepackage{amsthm}
\usepackage{amssymb}
\usepackage{stmaryrd}
\usepackage{setspace}
\usepackage{xargs}[2008/03/08]
\usepackage[unicode=true,pdfusetitle,
 bookmarks=true,bookmarksnumbered=false,bookmarksopen=false,
 breaklinks=false,pdfborder={0 0 1},backref=false,colorlinks=true]
 {hyperref}
\hypersetup{
 citecolor=blue,linkcolor=magenta,anchorcolor=green}

\makeatletter
\numberwithin{equation}{section}
\numberwithin{figure}{section}
\usepackage{changebar}
\providecommand{\lyxadded}[3]{}

\renewcommand{\lyxadded}[3]{
  {\protect\cbstart\color{lyxadded}{}#3\protect\cbend}
}

\@ifundefined{date}{}{\date{}}
\usepackage{upgreek}
\newrefformat{cor}{Corollary \ref{#1}}
\newrefformat{prop}{Proposition \ref{#1}}
\newrefformat{conj}{Conjecture \ref{#1}}
\newrefformat{nthm}{Theorem \nameref{#1}}
\newrefformat{def}{Definition \ref{#1}}

\makeatother

\begin{document}
\global\long\def\set#1#2{\left\{  #1\, |\, #2\right\}  }%

\global\long\def\cyc#1{\mathbb{Q}\!\left[\zeta_{#1}\right]}%

\global\long\def\mat#1#2#3#4{\left(\begin{array}{cc}
#1 & #2\\
#3 & #4
\end{array}\right)}%

\global\long\def\Mod#1#2#3{#1\equiv#2\, \left(\mathrm{mod}\, \, #3\right)}%

\global\long\def\inv{^{\,\textrm{-}1}}%

\global\long\def\pd#1{#1^{+}}%

\global\long\def\fix#1{\mathtt{Fix}\!\left(#1\right)}%

\global\long\def\map#1#2#3{#1\!:\!#2\!\rightarrow\!#3}%

\global\long\def\Map#1#2#3#4#5{\begin{split}#1:#2  &  \rightarrow#3\\
 #4  &  \mapsto#5 
\end{split}
 }%

\global\long\def\partlat{\boldsymbol{\sqcap}}%

\global\long\def\fact#1#2{#1\slash#2}%

\global\long\def\Gal#1{\mathtt{Gal}\!\left(#1\right)}%

\global\long\def\fixf#1{\mathbb{Q}\!\left(#1\right)}%

\global\long\def\gl#1#2{\mathsf{GL}_{#2}\!\left(#1\right)}%

\global\long\def\SL{\mathrm{SL}_{2}\!\left(\mathbb{Z}\right)}%

\global\long\def\zn#1{\left(\mathbb{Z}/\!#1\mathbb{Z}\right)^{\times}}%

\global\long\def\sn#1{\mathbb{S}_{#1}}%

\global\long\def\aut#1{\mathrm{Aut\mathit{\left(#1\right)}}}%

\global\long\def\FA#1{\vert#1\vert}%

\global\long\def\FB#1{\mathtt{Z}^{2}\!(#1)}%

\global\long\def\FC#1{#1^{{\scriptscriptstyle \flat}}}%

\global\long\def\FD#1{#1^{{\scriptscriptstyle \times}}}%

\global\long\def\FE#1{\mathtt{x}_{#1}}%

\global\long\def\FF#1#2{\mathrm{Fix}_{#1}\left(#2\right)}%

\global\long\def\FI#1{#1_{{\scriptscriptstyle \pm}}}%

\global\long\def\cl#1{\mathscr{C}\negthinspace\ell\!\left(#1\right)}%

\global\long\def\bl#1{\mathcal{B}\ell\!\left(#1\right)}%

\global\long\def\FJ#1#2#3{\uptheta\!\Bigl[\!{#1\atop #2}\!\Bigr]\!\left(#3\right)}%

\global\long\def\ex#1{\mathtt{e}^{2\mathtt{i}\pi#1}}%

\newcommandx\exi[3][usedefault, addprefix=\global, 1=, 2=]{\mathtt{e}^{\mathtt{{\scriptscriptstyle \textrm{#1}}}\frac{#2\pi\mathtt{i}}{#3}}}%

\global\long\def\tw#1#2#3{\mathfrak{#1}_{#3}^{{\scriptscriptstyle \left(\!#2\!\right)}}}%

\global\long\def\ver#1{\mathtt{Ver}_{#1}}%

\global\long\def\stab#1#2{#1_{#2}}%

\global\long\def\cft{\mathscr{C}}%

\global\long\def\gal#1{\upsigma_{#1}}%

\global\long\def\galpi#1{#1}%

\global\long\def\ann#1{\ker\mathfrak{#1}}%

\global\long\def\cent#1#2{\mathtt{C}_{#1}\!\left(#2\right)}%

\global\long\def\ch#1{\boldsymbol{\uprho}_{#1}}%

\global\long\def\chb#1{\xi_{#1}}%

\global\long\def\chg{\mathfrak{g}}%

\global\long\def\tgr{G}%

\global\long\def\qd#1{\mathtt{d}_{#1}}%

\global\long\def\cw#1{\mathtt{h}_{#1}}%

\global\long\def\tcl{\textrm{trivial class}}%

\global\long\def\ccl{\mathtt{C}}%

\global\long\def\zcl{\mathtt{z}}%

\global\long\def\om#1{\omega\!\left(#1\right)}%

\global\long\def\rami#1{\mathtt{e}_{#1}}%

\global\long\def\zm{\mathtt{Z}}%

\global\long\def\cs#1{\left\llbracket #1\right\rrbracket }%

\global\long\def\coc#1{\vartheta_{\mathfrak{#1}}}%

\global\long\def\sp#1#2{\bigl\langle#1,#2\bigr\rangle}%

\global\long\def\irr#1{\mathtt{Irr}\!\left(#1\right)}%

\global\long\def\stb#1{\mathtt{I}{}_{#1}}%

\global\long\def\sqf#1{#1^{\circ}}%

\global\long\def\CHG{\textrm{twister}}%

\global\long\def\fc{\textrm{FC set}}%

\global\long\def\to#1{\boldsymbol{\upsigma}\!\left(#1\right)}%

\global\long\def\v{\mathtt{{\scriptstyle 0}}}%

\global\long\def\wm#1{\mathcal{P}\!\left(#1\right)}%

\global\long\def\fm{\mathtt{N}}%

\global\long\def\du#1{#1^{{\scriptscriptstyle \perp}}}%

\global\long\def\dhg{\du{\chg}}%

\global\long\def\spr{\textrm{spread}}%

\global\long\def\gcl{\textrm{class}}%

\global\long\def\dcl{\textrm{block}}%

\global\long\def\vorb{\mathfrak{o}}%

\global\long\def\lat{\mathcal{\mathscr{L}}}%

\global\long\def\vera{\mathcal{V}}%

\global\long\def\svera#1{\vera_{#1}}%

\global\long\def\bd{\mathcal{D}_{\mathfrak{b}}}%

\global\long\def\rd#1{\boldsymbol{\updelta}_{#1}}%

\global\long\def\m#1{\mathfrak{\upmu}_{\mathfrak{#1}}}%

\global\long\def\zent#1{\mathtt{Z}\!\left(#1\right)}%

\global\long\def\voa{\mathcal{\mathbb{V}}}%

\global\long\def\bfm{\Delta}%

\global\long\def\cech#1#2{\boldsymbol{\varpi}_{#1}\!\left(#2\right)}%

\global\long\def\clchar#1#2{\updelta_{#1}\!\left(#2\right)}%

\global\long\def\vb{\mathfrak{b}_{{\scriptscriptstyle 0}}}%

\global\long\def\bch#1#2{\boldsymbol{\upchi}_{\mathfrak{#1}}\!\left(#2\right)}%

\global\long\def\prim{\mathtt{X}}%

\global\long\def\sig{\FA{\cft}}%

\global\long\def\beq#1{\boldsymbol{\Lambda}_{#1}}%

\global\long\def\twa#1{\boldsymbol{\nabla}\!{}_{#1}}%

\global\long\def\dg#1{\hat{#1}}%

\global\long\def\loc{\lat_{\mathtt{loc}}}%
\global\long\def\intlat{\lat_{\mathtt{int}}}%

\global\long\def\spp#1#2{\mathcal{P}_{#1}\!\left(#2\right)}%

\global\long\def\ifc{\boldsymbol{\Theta}}%
\global\long\def\sfc{\sqrt{\ifc}}%

\global\long\def\rcl{\mathtt{\mathsf{R}}}%
\global\long\def\ramcl{\textrm{Ramond~class}}%

\title{Character rings and fusion algebras}
\bigskip
\author{Peter Bantay}
\bigskip
\begin{abstract}
We present an overview of the close analogies between the character
rings of finite groups and the fusion rings of rational conformal
models, which follow from general principles related to orbifold deconstruction.
\end{abstract}

\curraddr{Institute for Theoretical Physics, E\" otv\" os L\' or\' and University,  H-1117 Budapest, P\' azm\' any P\' eter s. 1/A, Hungary}
\email{bantay@general.elte.hu}
\maketitle

\section{Introduction}

Analogies between group representation theory and 2D conformal field
theory have been noticed by several authors over the years. Some of
these have a natural interpretation because of the group theoretic
origin of the relevant conformal models (e.g. WZNW models based on
affine Lie algebras \cite{Borcherds1,Kac}, or holomorphic orbifolds
based on finite groups \cite{DV3}), but in other cases the relation
is less obvious. A new approach to the subject is provided by recent
advances in orbifold deconstruction, and the aim of the present note
is to give a sketchy overview of the relevant results.

Orbifold deconstruction \cite{Bantay2018a,Bantay2018} is a procedure
aimed at recognizing whether a given 2D conformal model is an orbifold
\cite{Dixon_orbifolds1,FLM1} of another one, and if so, to identify
(up to isomorphism) the relevant twist group and the original model.
The basic ideas have been described in \cite{Bantay2018a,Bantay2018},
focusing on the conceptual issues, while (part of) the relevant mathematical
details have been discussed in \cite{Bantay2019}. The basic observation
is that every orbifold has a distinguished set of primaries, the so-called
vacuum block, consisting of the descendants of the vacuum, and that
this vacuum block has quite special properties: it is closed under
the fusion product, and all its elements have integral conformal weight
and quantum dimension. Such sets of primaries, called twisters because
of their relation to twisted boundary conditions, provide the input
for the deconstruction procedure: each twister corresponds to a different
deconstruction, possibly with a different twist group and/or deconstructed
model.

Most of the basic notions related to twisters may be formulated in
the much more general setting of sets of primaries closed under the fusion product \cite{Bantay2019}, called \emph{$\fc$s} for short. As it turns out,
these behave in many ways as character rings of finite groups, especially
those -- the \emph{integral} ones --  all of whose elements have integral quantum dimension. In particular,
one may show that the collection of all $\fc$s of a conformal model
form a modular lattice (with the integral ones forming a sublattice),
allowing the generalization to $\fc$s of such group theoretic notions
as nilpotency, solubility, being Abelian, and so on. Moreover, there
is a well-defined notion of center and of central extensions, which
fit perfectly in the above mentioned analogy with group theory. Of
course, all this is completely natural for twisters, which are nothing
but character rings of twist groups according to the general principles
of orbifold deconstruction, but their meaning for general $\fc$s
still needs to be clarified. In the sequel, we shall try to sketch the highlights
of this circle of questions.

\section{FC sets and their classes\label{sec:FC}}

Let's consider a rational unitary conformal model \cite{DiFrancesco-Mathieu-Senechal,Ginsparg1988}:
we'll denote by $\qd p$ and $\mathtt{h}_{p}$ the quantum dimension
and conformal weight of a primary $p$, and by $\fm\!\left(p\right)$
the associated fusion matrix, whose matrix elements are given by the
fusion rules
\begin{equation}
\left[\fm\!\left(p\right)\right]_{q}^{r}=N_{pq}^{r}\label{eq:fmdef}
\end{equation}
The fusion matrices generate a commutative matrix algebra over $\mathbb{C}$,
the Verlinde algebra $\vera$, whose irreducible representations $\ch p$,
all of dimension $1$, are in one-to-one correspondence with the primaries
of the model.

An\emph{ $\mathcal{\fc}$ }is a set $\mathfrak{g}$ of primaries containing
the vacuum that is fusion closed, which means that 
\begin{equation}
\sum_{\gamma\in\chg}N_{\alpha\beta}^{\gamma}\qd{\gamma}=\qd{\alpha}\qd{\beta}\label{eq:fccrit}
\end{equation}
for all $\alpha,\beta\!\in\!\chg$. Taking into account that quantum
dimensions are positive numbers, this is tantamount to the requirement
that $N_{\alpha\beta}^{\gamma}\!>\!0$ and $\alpha,\beta\!\in\!\chg$
implies $\gamma\!\in\!\chg$.

The fusion matrices $\fm\!\left(\alpha\right)$ for $\alpha\!\in\!\mathfrak{g}$
generate a subalgebra $\svera{\mathfrak{g}}$ of the Verlinde algebra.
Because $\vera$ is commutative, the irreducible representations of
this subalgebra coincide with the different restrictions of the representations
$\ch p$ of $\vera$. A $\mathfrak{g}$-\emph{$\gcl$} is defined
to be the set $\ccl$ of all those primaries $p$ whose associated
representations $\ch p$ coincide when restricted to the subalgebra
$\svera{\mathfrak{g}}$; we shall denote by $\ch{\ccl}$ this common
restriction. Clearly, the collection $\cl{\mathfrak{g}}$ of $\mathfrak{g}$-$\gcl$es
(whose cardinality equals that of $\chg$) provides a partition of
the set of all primaries.

The first analogy with character rings of finite groups comes from
the existence of the orthogonality relations
\begin{equation}
\sum_{\ccl\in\cl{\mathfrak{g}}}\!\frac{\alpha\!\left(\ccl\right)\overline{\beta\!\left(\ccl\right)}}{\cs{\ccl}}=\begin{cases}
1 & \textrm{if }\alpha\!=\!\beta;\\
0 & \mathcal{\textrm{otherwise}}
\end{cases}\label{ortho1}
\end{equation}
 for $\alpha,\beta\!\in\!\mathfrak{g}$, and
\begin{equation}
\sum_{\alpha\in\chg}\alpha\!\left(\ccl_{1}\right)\overline{\alpha\!\left(\ccl_{2}\right)}=\begin{cases}
\cs{\ccl_{1}} & \textrm{if }\ccl_{1}\!=\!\ccl_{2};\\
0 & \mathcal{\textrm{otherwise}}
\end{cases}\label{eq:ortho2}
\end{equation}
for $\ccl_{1},\ccl_{2}\!\in\!\cl{\mathfrak{g}}$, where
\begin{equation}
\cs{\ccl}=\dfrac{\sum\nolimits _{p}\qd p^{2}}{\sum\limits _{p\in\ccl}\qd p^{2}}\label{eq:extdef}
\end{equation}
is the \emph{extent} of the class $\ccl\!\in\!\cl{\mathfrak{g}}$, and  $\alpha\!\left(\ccl\right)\!=\!\ch{\ccl}\!\left(\alpha\right)$.

The $\gcl$ containing the vacuum primary is of special importance:
we shall denote it by $\du{\chg}$ and call it the\emph{ $\tcl$}.
Note that $\alpha\!\left(\du{\chg}\right)\!=\!\qd{\alpha}$ for $\alpha\!\in\!\chg$,
and 
\begin{equation}
\cs{\du{\chg}}=\sum_{\alpha\in\chg}\qd{\alpha}^{2}\label{eq:spread}
\end{equation}
One can show that the trivial class maximizes both extent and the product
of size and extent, i.e. $\cs{\ccl}\!\leq\!\cs{\du{\mathfrak{g}}}$
and $\FA{\ccl}\!\cs{\ccl}\!\leq\!\FA{\du{\mathfrak{g}}}\!\cs{\du{\mathfrak{g}}}$
for every class $\ccl\!\in\!\cl{\mathfrak{g}}$.

A most important property of the trivial class is the \emph{product
rule}: if $N_{pq}^{r}\!>\!0$ for some $p\!\in\!\du{\chg}$, then
the primaries $q$ and $r$ belong to the same $\mathfrak{g}$-$\gcl$.
It follows at once that the trivial class $\du{\chg}$ is itself an
$\fc$, the \emph{dual} of $\mathfrak{g}$, hence all previous notions
and results go over verbatim to it. In particular, the set of all
primaries is partitioned into $\du{\mathfrak{g}}$-$\gcl$es, which
we shall call $\mathfrak{g}$-\emph{blocks} (or simply blocks) to
avoid confusion with $\mathfrak{g}$-$\gcl$es. It follows from what
has been said for classes that the set $\bl{\mathfrak{g}}$ of $\mathfrak{g}$-blocks
provides a partition of the set of all primaries. A simple argument
shows that the primaries $p$ and $q$ belong to the same $\mathfrak{g}$-$\dcl$
iff $N_{\alpha p}^{q}\!>\!0$ for some $\alpha\!\in\!\chg$, and in particular,
the $\mathfrak{g}$-block containing the vacuum is $\mathfrak{g}$
itself, that is $\du{\left(\du{\chg}\!\right)\!}\!=\!\chg$. This
illustrates the inherent duality of $\fc$s: $\chg$ and its dual
$\du{\chg}$ determine each other, while their extents are, roughly
speaking, reciprocal, since the product $\cs{\chg}\!\cs{\du{\chg}}$
can be shown to be the same for every $\fc$ $\mathfrak{g}$. This
duality implies that any result proven about classes gives a corresponding
result about blocks, and \emph{vice versa}, a seemingly trivial observation
that turns out to be quite useful.

The inclusion relation makes the collection $\lat$ of $\fc$s partially
ordered, with maximal element the $\fc$ containing all primaries,
and minimal element the trivial $\fc$ consisting of the vacuum primary
solely. Because the intersection of two $\fc$s is obviously an $\fc$
again, $\lat$ is actually a finite lattice, and one may show that
the join $\mathfrak{g}\!\vee\!\mathfrak{h}$ of the $\fc$s $\mathfrak{g}$
and $\mathfrak{h}$ (the smallest $\fc$ that contains both of them)
is the dual of the intersection of their duals, i.e.
\begin{equation}
\mathfrak{g}\!\vee\!\mathfrak{h}=\du{\left(\du{\mathfrak{g}}\cap\du{\mathfrak{h}}\right)}\label{eq:joindef}
\end{equation}
 If $\mathfrak{g}$ and $\mathfrak{h}$ are $\fc$s such that $\mathfrak{h\!\subseteq\!\mathfrak{g}}$,
then every $\mathfrak{h}$-class is a union of $\mathfrak{g}$-classes,
in particular $\du{\mathfrak{g}}\!\subseteq\!\du{\mathfrak{h}}$,
and every $\mathfrak{g}$-block is a union of $\mathfrak{h}$-blocks;
moreover, the number of $\mathfrak{g}$-classes contained in $\du{\mathfrak{h}}$
equals the number of $\mathfrak{h}$-blocks contained in $\mathfrak{g}$.
It follows that the map sending each $\fc$ to its dual is an isomorphism
between the lattice $\lat$ and its dual. An important consequence
of the above results, crucial from the viewpoint of orbifold deconstruction,
is that $\lat$ is a modular (even Arguesian) lattice, but usually
neither distributive nor complemented. A better understanding of the
lattice theoretic properties of $\lat$ would be highly desirable.

For an $\fc$ $\mathfrak{g}$, the of restriction of the indices of
the fusion matrices $\fm\!\left(\alpha\right)$ to the primaries belonging
to a given $\dcl$ $\mathfrak{\mathfrak{b}\!\in\!\bl{\mathfrak{g}}}$
results in non-negative integer matrices $\fm_{\mathfrak{b}}\!\left(\alpha\right)$
that form a representation $\bfm_{\mathfrak{b}}$ of the subalgebra
$\svera{\mathfrak{g}}$. This representation decomposes into a direct
sum of the irreducible representations $\ch{\ccl}$, and the \emph{overlap}
$\sp{\mathfrak{b}}{\ccl}$ of $\mathfrak{b}$ with the class $\ccl\!\in\!\cl{\mathfrak{g}}$
is defined as the multiplicity of $\ch{\ccl}$ in the irreducible
decomposition of $\bfm_{\mathfrak{b}}$. The overlap $\sp{\mathfrak{b}}{\ccl}$
may be shown to equal the rank of the minor of the modular $S$-matrix
obtained by restricting the row indices to $\mathfrak{b}\!\in\!\bl{\mathfrak{g}}$
and the column indices to $\ccl\!\in\!\cl{\mathfrak{g}}$. Alternatively,
one has the expression
\begin{equation}
\sp{\mathfrak{b}}{\ccl}=\sum_{p\in\mathfrak{b}}\sum_{q\in\ccl}\FA{S_{pq}}^{2}\label{eq:overlapdef}
\end{equation}
In particular, this means that $\sp{\mathfrak{b}}{\ccl}\!=\!0$ implies
$S_{pq}\!=\!0$ for all $p\!\in\!\mathfrak{b}$ and $q\!\in\!\ccl$,
explaining the appearance of large blocks of zeroes in the modular
$S$-matrix of many conformal models. Let's note that for twisters
(to be introduced in \prettyref{sec:Local-sets-and}) the overlap
has another, more profound group theoretic interpretation \cite{Bantay2018a,Bantay2018}.

\section{Central extensions\label{sec:The-center}}

As mentioned previously, the extent $\cs{\ccl}$ of a class $\ccl\!\in\!\cl{\mathfrak{g}}$
is bounded from above by the extent of the trivial class: $\cs{\ccl}\!\leq\!\cs{\du{\mathfrak{g}}}$.
Those classes  that saturate this bound -- the \emph{central} ones -- form the \emph{center} $\zent{\mathfrak{g}}\!=\!\set{\zcl\!\in\!\cl{\mathfrak{g}}}{\cs{\zcl}\!=\!\cs{\du{\chg}}}$
of the $\fc$ $\mathfrak{g}\!\in\!\lat$. Since $\du{\mathfrak{g}}\!\in\!\zent{\mathfrak{g}}$,
the center is never empty, and a straightforward argument shows that
$\zcl\!\in\!\zent{\mathfrak{g}}$ iff $\FA{\alpha\!\left(\zcl\right)\!}\!=\!\qd{\alpha}$
for all $\alpha\!\in\!\chg$. One may show that $\zent{\mathfrak{g}}\!=\!\mathfrak{g}$,
i.e. all classes are central precisely when $\qd{\alpha}\!=\!1$
for all $\alpha\!\in\!\mathfrak{g}$. Such $\fc$s are called (for
obvious reasons) \emph{Abelian;} in the language of 2D CFT, they correspond
to groups of simple currents.

The following generalization of the product rule holds: if the primary
$p$ belongs to the $\gcl$ $\ccl\!\in\!\cl{\mathfrak{g}}$ and $q$
belongs to the central $\gcl$ $\zcl\!\in\!\zent{\mathfrak{g}}$,
then all primaries $r$ for which $N_{pq}^{r}\!>\!0$ belong to the
same $\mathfrak{g}$-class, denoted $\zcl\ccl$. What is more, if
$\zcl_{1},\!\zcl_{2}\!\in\!\zent{\mathfrak{g}}$ are central classes,
then $\zcl_{1}\zcl_{2}$ is also central, and $\zcl_{1}\zcl_{2}\!=\!\zcl_{2}\zcl_{1}$,
hence the center $\zent{\mathfrak{g}}$ of an $\fc$ $\mathfrak{g}\!\in\!\lat$
is an Abelian group, and since $(\zcl_{1}\zcl_{2})\ccl\!=\!\zcl_{1}(\zcl_{2}\ccl)$
for any $\gcl$ $\ccl\!\in\!\cl{\mathfrak{g}}$, the center permutes
the $\mathfrak{g}$-classes.

\global\long\def\zquot#1#2{\mathfrak{#1}/\!#2}%
\global\long\def\mzq#1#2{#1^{[#2]}}%

For a subgroup $Z\!<\!\zent{\mathfrak{g}}$ of the center, $\zquot{\mathfrak{g}}Z\!=\!\set{\alpha\!\in\!\mathfrak{g}}{\alpha\!\left(\zcl\right)\!\!=\!\qd{\alpha}\textrm{ for all }\zcl\!\in\!Z}$
is again an $\fc$, the \emph{central quotient }of $\mathfrak{g}$
by $Z$. The usefulness of this notion rests on the explicit knowledge
of its structure, for one has complete control over its classes and
blocks in terms of those of $\mathfrak{g}$, in particular its dual
is given by
\begin{equation}
\du{\left(\zquot{\mathfrak{g}}Z\right)}=\bigcup_{\zcl\in Z}\zcl\label{eq:zentextdual}
\end{equation}
Moreover, if $\mathfrak{h}\!\in\!\lat$ is such that $\zquot{\mathfrak{g}}Z\!\!\subseteq\!\!\mathfrak{h}\!\subseteq\!\mathfrak{g}$,
then $\mathfrak{h}\!=\!\zquot{\mathfrak{g}}H$ for a suitable subgroup
$H\!<\!Z$. It follows that there is an order reversing one-to-one
correspondence between central quotients of $\mathfrak{g}\!\in\!\lat$
and subgroups of its center $\zent{\mathfrak{g}}$.

Given an $\fc$ $\mathfrak{g}$, it is natural to ask whether it is
a central quotient of another $\fc$. Given an Abelian group $A$,
an $A$-extension of $\mathfrak{g}$ is an $\fc$ $\mathfrak{h}\!\in\!\lat$
such that $\zquot{\mathfrak{h}}Z\!=\!\mathfrak{g}$ for some central
subgroup $Z\!<\!\zent{\mathfrak{h}}$ isomorphic to $A$. One may
show that the different $A$-extensions of $\mathfrak{g}$ are in
one-to-one correspondence with subgroups of $\zent{\du{\mathfrak{g}}}$
isomorphic to $A$, and in particular, every $\mathfrak{g}\!\in\!\lat$
has a maximal central extension, namely $\du{\left(\zquot{\du{\mathfrak{g}}}{\,\zent{\du{\mathfrak{g}}}}\right)}$,
the dual of the maximal central quotient of $\du{\mathfrak{g}}$.
Note that, while central quotients are related to groups of central
classes, central extensions have a similar relation to groups of central
blocks, i.e. blocks $\mathfrak{b}\!\in\!\bl{\mathfrak{g}}$ that satisfy
$\cs{\mathfrak{b}}\!=\!\cs{\mathfrak{g}}$, which form the center
$\zent{\du{\mathfrak{g}}}$ of the dual of $\mathfrak{g}$. An interesting
observation is that any class (hence any block) that contains a simple
current (a primary of dimension $1$) is automatically central, but
the converse need not be true.

\section{Local $\protect\fc$s\label{sec:Local-sets-and}}

It follows from the results discussed in \prettyref{sec:FC} that
every $\mathfrak{g}$-class is a union of $\mathfrak{g}$-blocks precisely
when $\mathfrak{g}\!\subseteq\!\du{\mathfrak{g}}$. It turns out that
such $\fc$s play a basic role in orbifold deconstruction \cite{Bantay2018a,Bantay2018},
hence they deserve a special name: we'll call them \emph{local} $\fc$s.
The point is that the vacuum block of an orbifold model (the set of
primaries that originate from the vacuum) is an $\fc$ whose classes
correspond to the different twisted sectors, i.e. collections of twisted
modules whose twist element belong to the same conjugacy class, while
its blocks correspond to orbits of twisted modules. Since the conjugacy
class of a twist element is the same for all twisted modules on the
same orbit, every block should be included in a well-defined class,
hence the vacuum block, considered as an $\fc$, should be local by
the above. Note that, while the intersection of local $\fc$s is clearly
local, this is not necessarily the case for their join, hence they
do not form a lattice.

Actually, the vacuum block of an orbifold belongs to a special class
of local $\fc$s, termed \emph{twisters} because of their relation
to twisted boundary conditions, characterized by all of their elements
having integral conformal weight. Indeed, since elements of the vacuum
block descend from the vacuum primary, they all have integral conformal
weight, hence the vacuum block is necessarily a $\CHG$.

One may show that a local $\fc$ is either itself a $\CHG$, or a
$\mathbb{Z}_{2}$-extension of a $\CHG$. More precisely, every local
$\fc$ that is not a $\CHG$ has a distinguished central class $\rcl$,
the so-called \emph{Ramond class}, such that the corresponding central quotient
is a $\CHG$. The rationale for this appellation is that, in a suitable
fermionic generalization of orbifold deconstruction, the blocks contained
in the trivial class account for the Neveu-Schwarz (bosonic) sector
of the deconstructed model, while those contained in the $\ramcl$
describe the fermionic (Ramond) sector. This is substantiated by the
observation that a $\dcl$ is contained in the $\ramcl$ precisely
when the conformal weights of its elements differ by integers, and
that the number of blocks contained in the trivial and $\ramcl$es
(the number of bosonic and fermionic degrees of freedom) are equal.

Local $\fc$s have many striking properties. For example, one may
show that all of their elements have (rational) integer dimension,
and either integer or half-integer conformal weight\footnote{Note that the converse is not true: there are many $\fc$s in which
all conformal weights belong to $\frac{1}{2}\mathbb{Z}$, but are
nevertheless not local; on the other hand, the integrality of conformal
weights implies locality, hence every twister is automatically local.}. Ultimately, all this follows from the observation that, if the primaries
$\alpha,\beta$ belong to a local $\fc$ and $N_{\alpha\beta}^{\gamma}\!>\!0$,
then the conformal weight $\cw{\gamma}$ differs by an integer from
the sum $\cw{\alpha}\!+\!\cw{\beta}$. From a categorical point of
view this means that the elements of a local $\fc$ are the simple
objects of a Tannakian subcategory of the modular tensor category
associated to the model. According to a result of Deligne \cite{Deligne1990},
such a category is (tensor-)equivalent to the representation category
of some finite group, hence the associated subalgebra may be identified
with the character ring of that group. This has a natural interpretation
in terms of orbifold deconstruction: the elements of the vacuum block
correspond to irreducible representations of the twist group, hence
their fusion rules describe the decomposition of tensor products of
the latter, implying that the subalgebra $\svera{\mathfrak{g}}$ of
the Verlinde algebra is nothing but the character ring of the twist
group.

It follows from the above considerations that results from character
theory \cite{Isaacs,Lux-Pahlings,Serre} should go over to local $\fc$s,
and this observation allows the generalization of many group theoretic
concepts to arbitrary $\fc$s, providing a host of non-trivial conjectural
results that seem to hold in greater generality. In this vein one
can generalize \cite{Bantay2019} to $\fc$s such concepts as nilpotency,
(super)solubility, and so on. For example, an $\fc$ $\mathfrak{g}$
is \emph{nilpotent} if it can be obtained from the trivial $\fc$
by a sequence of central extensions. The rationale for this terminology
is that if $\mathfrak{g}$ is local, hence the algebra $\svera{\mathfrak{g}}$
is isomorphic to the character ring of some finite group $G$, then
$\mathfrak{g}$ is nilpotent precisely when the corresponding group
$G$ is. One may show that if the $\fc$ $\mathfrak{g}$ is nilpotent,
then $\cs{\du{\mathfrak{g}}}\!\in\!\mathbb{Z}$, and this in turn
implies that the quantum dimension $\qd{\alpha}$ of any element 
$\mathfrak{\alpha\!\in\!g}$ is either an integer or the square root
of an integer\footnote{That the latter possibility can occur is exemplified by the maximal
$\fc$ of the Ising model (the minimal Virasoro model of central charge
$\tfrac{1}{2}$), which is nilpotent while having a primary of dimension
$\sqrt{2}$.}. We conjecture that many results about (finite) nilpotent groups
would carry over to nilpotent $\fc$s, e.g.  if $\mathfrak{g}$ is
nilpotent and $d$ is an integer dividing $\cs{\du{\mathfrak{g}}}$,
then there would exist an $\fc$ $\mathfrak{h}\!\subseteq\!\mathfrak{g}$
such that $\cs{\du{\mathfrak{h}}}\!=\!d$.

As explained above, for a local $\fc$ $\mathfrak{g}$ the algebra
$\svera{\mathfrak{g}}$ is isomorphic to the character ring of some
finite group, hence usual properties of character rings \cite{Isaacs,Lux-Pahlings,Serre}
should apply to it. This means in particular that
\begin{enumerate}
\item the extent $\cs{\ccl}$ of any $\mathfrak{g}$-class $\ccl\!\in\!\cl{\mathfrak{g}}$
is a rational integer dividing $\cs{\du{\mathfrak{g}}}$;
\item the dimension $\qd{\alpha}$ of every $\alpha\!\in\!\mathfrak{g}$
is a rational integer dividing ${\displaystyle \frac{\cs{\du{\mathfrak{g}}}}{\FA{\zent{\mathfrak{g}}\!}}}$
;
\item if $\alpha\!\in\!\mathfrak{g}$ has dimension $\qd{\alpha}\!>\!1$,
then $\alpha\!\left(\ccl\right)\!=\!0$ for some class $\ccl\!\in\!\cl{\mathfrak{g}}$;
\item if the class $\ccl\!\in\!\cl{\mathfrak{g}}$ is such that ${\displaystyle \frac{\cs{\du{\mathfrak{g}}}}{\cs{\ccl}}}$
is coprime to $\qd{\alpha}$ for some $\alpha\!\in\!\mathfrak{g}$,
then either $\FA{\alpha\!\left(\ccl\right)\negmedspace}\!=\!\qd{\alpha}$
or $\alpha\!\left(\ccl\right)\!=\!0$.
\end{enumerate}
All these assertions are well-known properties of character rings,
e.g. the first one just states that the size of a conjugacy class
is an integer dividing the order of the group, while the second one
is Ito's celebrated theorem \cite{Isaacs}. What is really amazing
is that, as suggested by extensive computational evidence, these properties
(and many similar ones) seem to hold for a much larger class of $\fc$s
\cite{Bantay2019}, the so-called integral ones characterized by the
property that all of their elements have  integer dimension\footnote{Actually, they seem to hold (after suitable amendments) in the more
general case of $\fc$s whose elements have integer squared dimension.}. This is truly surprising, as one can show explicit examples of integral
$\fc$s which are not the character ring of any finite group. From
this point of view, one may consider integral $\fc$s as 'character
rings' of some suitable generalization of the group concept.

\section{Summary}

As we have seen, fusion closed sets of primaries of a conformal model
(or modular tensor category) have a fairly deep structure, generalizing
many aspects of the character theory of finite groups. Of course,
this is no accident, since the vacuum block of an orbifold model,
which corresponds on general grounds to the character ring of the
twist group, is a special variety of $\fc$, a so-called twister.
But it turns out that the parallel with character theory goes much
further, even for $\fc$s that have no group theoretic interpretation.
Many notions from group theory (like nilpotency, solubility, etc.)
may be generalized to arbitrary $\fc$s, and the corresponding properties
seem to go over almost verbatim in this more general setting. But
it should be stressed that $\fc$s are more than just some fancy generalization
of the group concept, since they possess genuinely new structures
and properties, which await careful study.

\bibliographystyle{plain}


\end{document}